\documentclass[11pt]{amsart}
\usepackage{amsmath,graphicx}
\newtheorem{thm}{Theorem}[section]
\newtheorem{lem}[thm]{Lemma}
\newtheorem{prop}[thm]{Proposition}
\newtheorem{cor}[thm]{Corollary}
\theoremstyle{remark}
\theoremstyle{definition}

\theoremstyle{plain}

\newcommand{\eps}{\varepsilon}
\newcommand{\Z}{{\mathbb{Z}}}

\newcommand{\R}{{\mathbb{R}}}

\newcommand{\T}{{\mathbb{T}}}

\newcommand{\FQ}{{\mathcal{F}_Q}}
\newcommand{\FI}{\FQ (I)}
\newcommand{\EE}{{\mathcal{E}}}
\newcommand{\II}{{\mathcal{I}}}
\newcommand{\JJ}{{\mathcal{J}}}

\newcommand{\dist}{\operatorname{dist}}

\numberwithin{equation}{section}

\begin{document}

\title[The statistics of the trajectory of a billiard]{The statistics of the trajectory of a certain billiard
in a flat two-torus}
\author[F.P. Boca, R.N. Gologan, A. Zaharescu]{Florin P. Boca, Radu N. Gologan and Alexandru Zaharescu}

\address{FPB and AZ: Department of Mathematics, University of
Illinois at \break Urbana-Champaign, Urbana, IL 61801, USA;
E-mail: fboca@math.uiuc.edu;\break zaharesc@math.uiuc.edu}

\address{Institute of Mathematics of the Romanian Academy,
P.O. Box 1-764, \break Bucharest RO-014700, Romania}

\address{RNG: Institute of Mathematics of the Romanian Academy,
P.O. Box 1-764, Bucharest RO-014700, Romania; E-mail:
Radu.Gologan@imar.ro}

\thanks{Research partially supported by ANSTI grant C6189/2000}

\begin{abstract}
We consider a billiard in the punctured torus obtained by removing
a small disk of radius $\eps >0$ from the flat torus $\T^2$,
with trajectory starting from the center of the puncture.
In this case the phase space is given by the range of the velocity
$\omega$ only. Let $\tilde{\tau}_\eps (\omega)$, and respectively
$\tilde{R}_\eps (\omega)$, denote the first exit time (length of
the trajectory), and respectively the number of collisions with the
side cushions when $\T^2$ is being identified with $[0,1)^2$.
We prove that the probability measures on $[0,\infty)$ associated
with the random variables $\eps \tilde{\tau}_\eps$ and
$\eps \tilde{R}_\eps$ are weakly convergent as $\eps \rightarrow 0^+$
and explicitly compute the densities of the limits.
\end{abstract}

\maketitle

\section{Introduction and main results}
Various ergodic and statistical properties of the periodic Lorentz
gas were studied during the last decades (see \cite{Ble}, \cite{BGW},
\cite{BuSi1}, \cite{BuSi2}, \cite{BSC1}, \cite{BSC2},
\cite{CG}, \cite{Ch1}, \cite{Ch2}, \cite{Dah}, \cite{DDG1},
\cite{DDG2}, \cite{FOK}, \cite{Ku}, \cite{Sin}, \cite{Sin90}
for a non-exhaustive list of references). In the case of
periodically distributed circular scatterers of radius $\eps \in
(0,\frac{1}{2})$ in $\R^2$, one considers the region
\begin{equation*}
Z_\eps =\{ x\in \R^2 \, ;\, \dist (x,\Z^2) \geq \eps \}
\end{equation*}
and the {\it first exit time} (called {\it free path length}
by some authors) defined as
\begin{equation*}
\tau_\eps (x,\omega)=\inf \{ \tau >0 \, ;\, x+\tau \omega \in
\partial Z_\eps \} ,
\end{equation*}
where $x\in Y_\eps =Z_\eps/\Z^2$ and $\omega$ belongs to
the unit circle $\T$, which will be steadily identified with
$[0,2\pi)$ throughout the paper.
Equivalently, one can consider the free motion of a point-like
particle in the billiard table $Y_\eps$ obtained by removing
pockets of the form of quarters of a circle of radius $\eps$. If
we identify $\T^2=\R^2/\Z^2$ with $[0,1)^2$, then $Y_\eps$ can be
regarded as a punctured two-torus. The reflection in the side
cushions of the table is specular and the trajectory between two
such reflections is rectilinear. Assume that the particle has
constant speed, say equal to $1$, and leaves the table when it
reaches one of the four pockets. In this setting $\tau_\eps
(x,\omega)$ denotes the exit time from the table (or equivalently
the length of the trajectory).

\begin{figure}[ht]
\begin{center}
\includegraphics*[scale=0.7, bb=0 0 190 190]{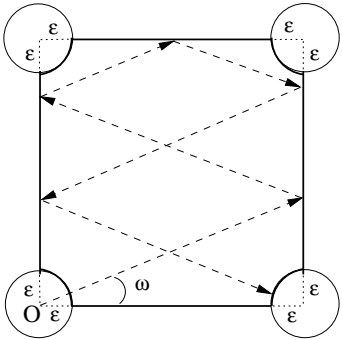}
\end{center}
\caption{The trajectory of the billiard} \label{Figure1}
\end{figure}

In this paper we will study the case where the trajectory
starts at the origin $O=(0,0)$ with initial velocity $\omega$.
In this case the
first exit time is $\tilde{\tau}_\eps (\omega) :=\tau (O,\omega)$
and one averages over $\omega$ only. We shall give
very precise estimates on the average of $\tilde{\tau}_\eps$ and
of the number $\tilde{R}_\eps (\omega)$ of collisions of the
particle with the side cushions as $\eps \rightarrow 0^+$. The
related problem of estimating the moments
\begin{equation*}
c_r=\frac{1}{2\pi} \int\limits_0^{2\pi} \tilde{\tau}_\eps^r
(\omega)\, d\omega ,\qquad r>0,
\end{equation*}
was raised by Ya.G. Sinai in 1981. An answer was given in
\cite{BGZ}, where it was proved for any interval $I\subseteq
[0,\frac{\pi}{4}]$ and any (small) $\delta>0$ that
\begin{equation}\label{1.1}
\eps^r \int\limits_I \tilde{\tau}_\eps^r (\omega)\, d\omega =
c_r \vert I\vert +O_{r,\delta}
(\eps^{\frac{1}{8}-\delta}) \quad \mbox{as}\ \eps \ \rightarrow
0^+,
\end{equation}
where
\begin{equation*}
c_r =\frac{2}{\zeta(2)} \int\limits_0^{1/2} \left( x\big(
x^{r-1}+(1-x)^{r-1}\big) +\frac{1-(1-x)^r}{rx(1-x)} -
\frac{1-(1-x)^{r+1}}{(r+1)x(1-x)} \right) dx .
\end{equation*}

Since all the probability measures $\tilde{\mu}_\eps$ defined by
\begin{equation*}
\tilde{\mu}_\eps (f)=\frac{1}{2\pi} \int\limits_0^{2\pi} f\big(
\eps \tilde{\tau}_\eps (\omega)\big)\, d\omega ,\qquad f\in C_c
\big( [0,\infty )\big) ,
\end{equation*}
have the support included into a common compact (see Lemma 3.1 in
\cite{BGZ} or Lemma 2.1 in this paper), the equality \eqref{1.1}
implies that $\tilde{\mu}_\eps \rightarrow \tilde{\mu}$ weakly
as $\eps \rightarrow 0^+$, and that the moments of
$\tilde{\mu}$ are given by
\begin{equation*}
\frac{1}{2\pi} \int\limits_0^{2\pi} \omega^r d\tilde{\mu} (\omega)
=c_r .
\end{equation*}
However, this does not directly provide an explicit formula
for the density of
$\tilde{\mu}$. The primary aim of this paper is to give a direct
proof for this convergence, computing in the meantime the density
of $\tilde{\mu}$ in closed form. For obvious symmetry reasons we
may only consider the interval $[0,\frac{\pi}{4}]$ instead of
$[0,2\pi]$. To state the main result we consider the repartition
function
\begin{equation*}
H_\eps (t)=\frac{4}{\pi} \left| \left\{ \omega \in \Big[
0,\frac{\pi}{4} \Big] \, ;\, \eps \tilde{\tau}_\eps (\omega) >t
\right\} \right|
\end{equation*}
of $\tilde{\mu}_\eps$ and the function
\begin{equation}\label{1.2}
\psi (x)=\frac{1-x}{x} \left( 1+\ln \frac{x}{1-x} \right),\qquad
x\in \Big[ \frac{1}{2},1\Big) .
\end{equation}
It is seen that $\psi$ is non-decreasing and satisfies
$\psi (\frac{1}{2})=1$, $\psi (1^-)=0$, and especially
\begin{equation}\label{1.3}
\int\limits_{1/2}^1 \psi (x)\, dx=\frac{\zeta(2)-1}{2} =
\frac{\pi^2}{12}-\frac{1}{2} \, .
\end{equation}

We can now state the main result.

\medskip

\begin{thm}\label{T1.1}
For each $t>0$, there exists $H(t)=\lim\limits_{\eps \rightarrow
0^+} H_\eps (t)$. Moreover, one has
\begin{equation*}
H_\eps (t)=H(t)+O_\delta (\eps^{\frac{1}{8}-\delta})
\end{equation*}
for any {\em (}small{\em )} $\delta >0$, where
\begin{equation*}
H(t)=\begin{cases} \vspace{0.2cm} \displaystyle
1-\frac{2t}{\zeta(2)} & \mbox{\rm
if $\  0<t<\frac{1}{2}$;} \\ \vspace{0.2cm}
\displaystyle \frac{2}{\zeta(2)}\int_t^1 \psi (x)\, dx &
\mbox{\rm if $\ \frac{1}{2}<t<1$;}
\\ 0 & \mbox{\rm if $\ t>1$.}
\end{cases}
\end{equation*}
In particular $\tilde{\mu}_\eps \rightarrow \tilde{\mu}$ weakly
as $\eps \rightarrow 0^+$, where $d\tilde{\mu}(t)=h(t) dt$ and
\begin{equation*}
h(t)=-H^\prime (t)=\frac{2}{\zeta(2)} \cdot \begin{cases} 1 &
\mbox{\rm if $\ 0<t\leq \frac{1}{2}$;} \\
\psi (t) & \mbox{\rm if $\  \frac{1}{2} \leq t\leq 1$;} \\ 0
& \mbox{\rm if $\ t\geq 1$.}
\end{cases}
\end{equation*}
\end{thm}

\begin{figure}[ht]
\includegraphics*[scale=0.72, bb=90 0 330 150]{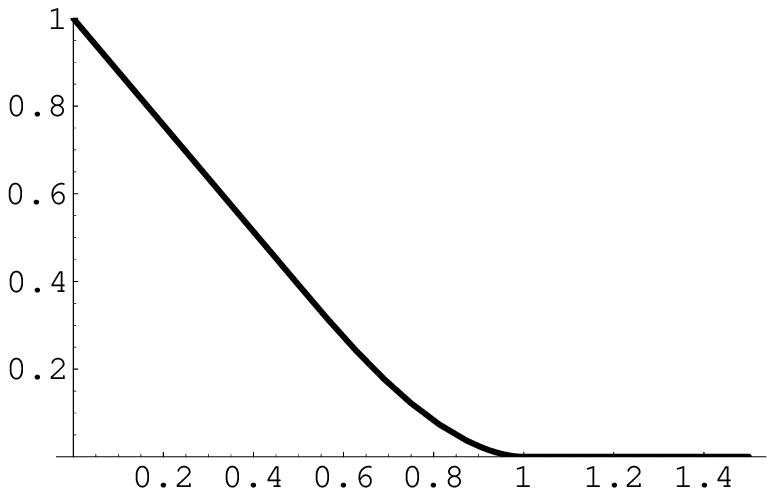}
\includegraphics*[scale=0.72, bb=90 0 330 150]{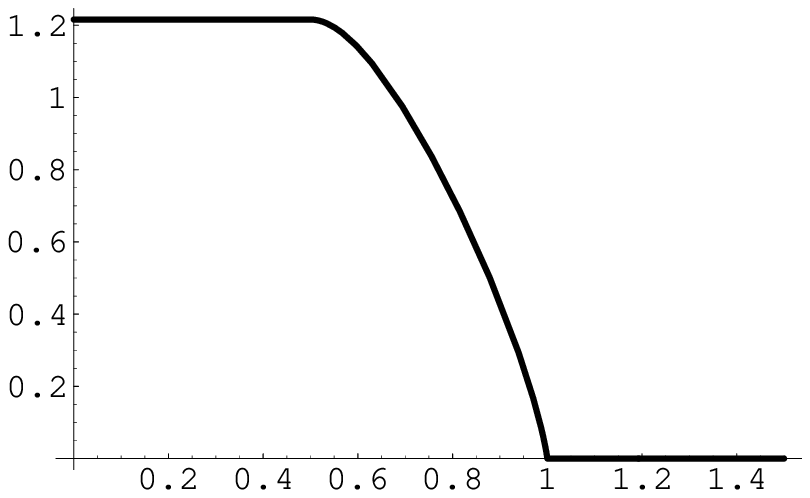}
\caption{The graphs of $H(t)$, and respectively of
$h(t)=-H^\prime (t)$}
\label{Figure2}
\end{figure}

This result deserves some comments. Firstly, since $\tilde{\mu}$
is a probability measure on $[0,1]$, one must have that
$\int_0^1 h(t)\, dt=1$. But this amounts to
\begin{equation*}
\frac{1}{2}+\int\limits_{1/2}^1 \psi(t)\, dt=\frac{\zeta(2)}{2}\, ,
\end{equation*}
that is to formula \eqref{1.3}, and we got a first proof of that formula.

However, one can prove \eqref{1.3} in other ways. For instance, writing
$\psi(x)=\frac{1}{x}\ln \frac{1}{1-x}+\ln (1-x)+\frac{1+\ln x}{x}
-1-\ln x$, integrating each term from $t$ to $1$ and using
the Taylor expansion for $\frac{1}{x}\ln \frac{1}{1-x}$, we gather
\begin{equation}\label{1.4}
\int\limits_t^1 \psi(x)\, dx={\mathrm {Li}}_2 (1)
-{\mathrm {Li}}_2(t)-\frac{\ln^2 t}{2}
+(1-t)\left( \ln \frac{1-t}{t} -1\right) ,
\end{equation}
where ${\mathrm{Li}}_2(t)=\sum_1^\infty \frac{t^n}{n^2}$
denotes the dilogarithm
function. Plugging now $t=\frac{1}{2}$ in \eqref{1.4} and using the
formula (cf. relation (1.6) in \cite{Lew})
\begin{equation}\label{1.5}
{\mathrm{Li}}_2 \left( \frac{1}{2}\right)
=\frac{\pi^2}{12}-\frac{\ln^2 2}{2}=
\frac{\zeta(2)-\ln^2 2}{2}
\end{equation}
we arrive again at relation \eqref{1.3}.
Thus, the fact that $\tilde{\mu}$
has mass one is equivalent to the (non-trivial) equality
\eqref{1.5} concerning the dilogarithm.

Secondly, Theorem \ref{T1.1} should be compared with the most recent
results of the first two authors (\cite{BZ}), who proved the existence
of the limit repartition
\begin{equation*}
\frac{4}{\pi} \left| \left\{ (x,\omega)\in [0,1)^2 \times \Big[
0,\frac{\pi}{4}\Big] \, ;\, \eps \tau_\eps (x,\omega)>t\right\}
\right| \qquad \mbox{\rm as $\ \eps \rightarrow 0^+$.}
\end{equation*}
This repartition function is no longer
compactly supported; it was proved in \cite{BZ} to be of the form
\begin{equation*}
4\zeta(2)\sum\limits_{n=1}^\infty \frac{2^n-1}{n^2(n+1)^2(n+2)t^n}
\ +O_\delta (\eps^{\frac{1}{8}-\delta}),\qquad \forall \, t\geq 2.
\end{equation*}
Interestingly, the density $h$ coincides\footnote{after scaling
$t$ by $2$ since the factor $\eps$ is replaced by $2\eps$
in \cite{BZ}} for $t\in [0,1]$ with the density of the limit
measure of the geometric free path, proved to exist in the
small-scatterer limit in \cite{BZ}, and computed independently
in \cite{Dah} and \cite{BZ}.
Another interesting feature of the free path lengths in the
small-scatterer limit is that their limit repartition functions
are closely related with dilogarithms in the homogeneous
case treated in this paper, and with dilogarithms and possibly
trilogarithms in the nonhomogeneous case treated in \cite{BZ}, where
one averages over a phase space defined by both velocity and position,
and where the limit measure has a tail at $+\infty$.

We finally estimate the repartition function
\begin{equation*}
F_\eps (t)=\frac{4}{\pi} \left| \left\{ \omega \in \Big[
0,\frac{\pi}{4} \Big] \, ;\, \eps \tilde{R}_\eps (\omega)>t
\right\} \right|
\end{equation*}
of the random variable $\eps \tilde{R}_\eps$ in the
small-scatterer limit, and prove

\medskip

\begin{thm}\label{T1.2}
For each $t>0$, there exists $F(t)=\lim\limits_{\eps \rightarrow 0^+}
F_\eps (t)$. Moreover, one has
\begin{equation*}
F_\eps (t)=F(t)+O_\delta (\eps^{\frac{1}{8}-\delta})
\end{equation*}
for any {\em (}small{\em )} $\delta >0$, where
\begin{equation*}
F(t)=\frac{4}{\pi} \int\limits_0^{\pi/4} H \left(
\frac{t}{\cos \omega+\sin \omega} \right) d\omega
=\frac{4}{\pi} \int\limits_{\pi/4}^{\pi/2}
H\left( \frac{t}{\sqrt{2} \sin \omega} \right) d\omega.
\end{equation*}
\end{thm}

\medskip

In particular this implies that the probability measures
$\tilde{\nu}_\eps$ defined by
\begin{equation*}
\tilde{\nu}_\eps (f)=\frac{1}{2\pi} \int\limits_0^{2\pi}
f\big( \eps \tilde{R}_\eps (\omega)\big) \, d\omega ,
\qquad f\in C_c \big( [0,\infty)\big),
\end{equation*}
converge weakly to the probability measure $\tilde{\nu}$ with
repartition function $F$. It is clear that the support of
$\tilde{\nu}$ coincides with the interval $[0,\sqrt{2}]$.
Moreover, the function $F$ is linear on $[0,\frac{1}{2}]$
and on this interval one has
\begin{equation*}
\begin{split}
F(t) & =\frac{4}{\pi} \int\limits_{\pi/4}^{\pi/2}
\left( 1-\frac{t\sqrt{2}}{\zeta(2)\sin \omega} \right) d\omega
=1+\frac{4t\sqrt{2}}{\pi \zeta(2)}\cdot \ln \tan \frac{\pi}{8} \\
& =1-\frac{4\sqrt{2} \ln (1+\sqrt{2}) t}{\pi \zeta(2)} \, .
\end{split}
\end{equation*}

\medskip

One may replace the range $[0,\frac{\pi}{4}]$ of $\omega$ by
an arbitrary interval $I\subseteq [0,\frac{\pi}{4}]$, proving
the existence of the weak limit of the probability measures associated
with the random variables $\eps \tilde{\tau}_\eps$ and
$\eps \tilde{R}_\eps$ when $\omega \in I$, and also computing their
limits. This can be done only with minimal modifications in Section 4.

\bigskip

\section{Formulas for sectors ending at Farey points}
For each integer $Q\geq 1$, let $\FQ$ denote the set of Farey
fractions of order $Q$, that is the set of irreducible rational
numbers in $[0,1]$ with denominators less or equal
than $Q$. It is well-known
that if $\frac{a}{q}<\frac{a^\prime}{q^\prime}$ are consecutive
elements in $\FQ$, then one has
\begin{equation*}
a^\prime q-aq^\prime =1 \qquad \mbox{\rm and} \qquad q+q^\prime
>Q.
\end{equation*}
Conversely, if $q,q^\prime \in \{ 1,\dots,Q\}$ and $q+q^\prime
>Q$, then there exist $a\in \{ 1,\dots,q-1\}$ and $a^\prime \in \{
1,\dots,q^\prime -1 \}$ such that
$\frac{a}{q}<\frac{a^\prime}{q^\prime}$ are consecutive in $\FQ$.

For each interval $I\subseteq [0,1]$, we consider the set
\begin{equation*}
\FI =I\cap \FQ
\end{equation*}
of Farey fractions of order $Q$ from $I$ of cardinality
\begin{equation*}
N_{Q,I}=\# \FI =\frac{Q^2 \vert I\vert}{2\zeta(2)}+O(Q\ln Q) .
\end{equation*}

For each $h>0$ we consider the vertical scatterers of height $2h$,
\begin{equation*}
\tilde{V}_{q,a,h} =\{ q\} \times [a-h,a+h] ,\qquad (q,a)\in
\Z^{2\ast} ,
\end{equation*}
and the set
\begin{equation*}
\tilde{Z}_h =\bigcup\limits_{(q,a)\in \Z^{2\ast}} \hspace{-4pt}
\tilde{V}_{q,a,h} .
\end{equation*}
The {\it free path length} in this model is then given by
\begin{equation*}
l_h (\omega) =\inf \{ \tau >0\, ;\, (\tau \cos \omega,\tau \sin
\omega ) \in \tilde{Z}_h \} .
\end{equation*}
We denote by $R_h (\omega)$ the number of reflections in the side
cushions in the billiard model in the case of vertical scatterers
considered above.

We also define for each $\eps >0$, each interval $I\subseteq
[0,1]$ and each integer $Q\geq 1$ the quantities
\begin{equation*}
\tilde{H}_{\eps,I,Q}(t)=\bigg| \bigg\{ \omega \, ;\, \tan \omega
\in I, \ l_{\frac{1}{Q}} (\omega) >\frac{t}{\eps} \bigg\} \bigg|
\end{equation*}
and
\begin{equation*}
\tilde{F}_{\eps,I,Q}(t)=\left| \left\{ \omega \, ;\, \tan \omega
\in I,\ R_{\frac{1}{Q}} (\omega)>\frac{t}{\eps} \right\} \right| .
\end{equation*}

The repartition of Farey fractions in
$\FI$ will play a central role in the next section while proving
an asymptotic formula for $\tilde{H}_{\eps,I,Q}(t)$ as $\eps \rightarrow
0^+$, $\vert I\vert \rightarrow 0$ and
$Q\rightarrow \infty$ in a suitable way.

In the remainder of this section $Q\geq 1$ will be a fixed
integer. A key remark is that every line of slope between $0$ and
$1$ through the origin $O$ will necessarily intersect the set
\begin{equation*}
\mathcal{V}_Q =\bigcup\limits_{a/q \in \FQ}
\tilde{V}_{q,a,\frac{1}{Q}}
\end{equation*}
consisting in $N_Q=N_{Q,[0,1]}$ vertical scatterers of height
$\frac{2}{Q}$ centered at the integer points $(q,a)$ with
$\frac{a}{q}$ irreducible fraction in $\FQ$. This is contained in
the next statement which is essentially Lemma 3.1 in \cite{BGZ}.

\medskip

\begin{lem}\label{L2.1}
For any $\omega \in [0,\frac{\pi}{4}]$ one has
\begin{equation*}
\{ (\tau \cos \omega ,\tau \sin \omega)  \, ;\, \tau >0\} \cap
\mathcal{V}_Q \neq \emptyset .
\end{equation*}
\end{lem}

{\sl Proof.} Let $t_P$ denote the slope of the line $OP$. We use
the inequalities $q+q^\prime \geq Q+1 >\max (q,q^\prime)$ to infer
\begin{equation*}
t_A=\frac{a}{q}\leq t_{S^\prime} =\frac{a^\prime
-\frac{1}{Q}}{q^\prime} <t_N=\frac{a+\frac{1}{Q}}{q} \leq
t_{A^\prime} =\frac{a^\prime}{q^\prime} \, ,
\end{equation*}
where we set $A=(q,a)$, $A^\prime=(q^\prime,a^\prime)$,
$N=(q,a+\frac{1}{Q})$, $N^\prime
=(q^\prime,a^\prime+\frac{1}{Q})$, $S=(q,a-\frac{1}{Q})$,
$S^\prime =(q^\prime ,a^\prime -\frac{1}{Q})$,
$W=(q-\frac{1}{Q},a)$, $W^\prime =(q^\prime
-\frac{1}{Q},a^\prime)$. This clearly shows that for any $\omega
\in [0,\frac{\pi}{4}]$, the line of slope $\omega$ through the
origin will necessarily intersect $\mathcal{V}_Q$. \qed
\begin{figure}[ht]
\includegraphics*[scale=0.8, bb=0 0 230 145]{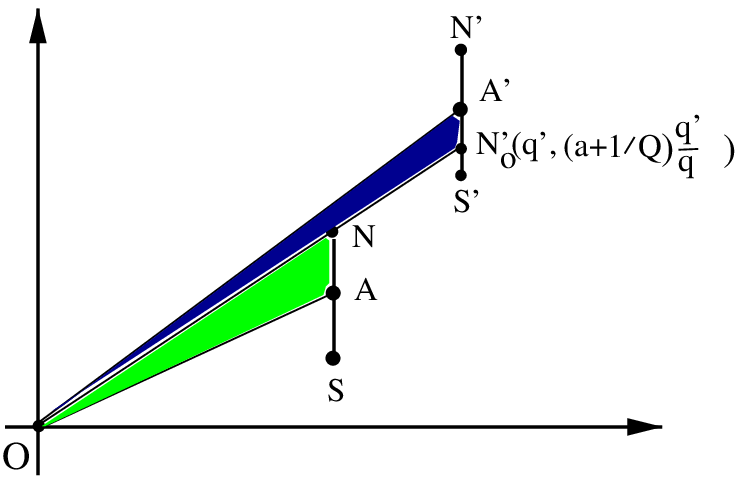}
\includegraphics*[scale=0.72, bb=0 0 230 180]{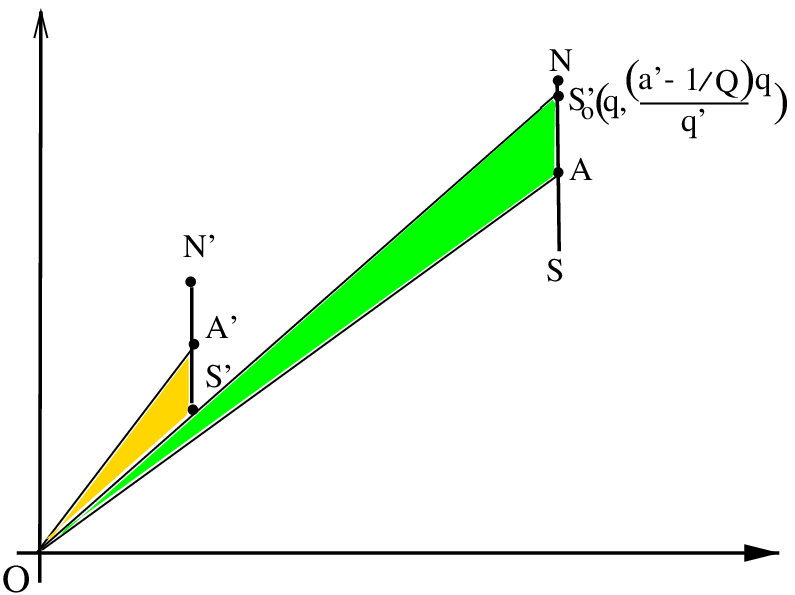}
\caption{The cases $q<q^\prime$ and respectively $q>q^\prime$}
\label{Figure3}
\end{figure}

\bigskip

We need a few more elementary things concerning Farey fractions.
Suppose next that
$\frac{a^\prime}{q^\prime}<\frac{a}{q}<\frac{a^{\prime
\prime}}{q^{\prime \prime}}$ are three consecutive fractions in
$\FQ$. Then the relation $aq^\prime -a^\prime q=1$ yields that
$q^\prime =\bar{a} \hspace{-3pt} \pmod{q}$, where $\bar{a}$
denotes the multiplicative inverse of $a \hspace{-3pt} \pmod{q}$.
Since $q^\prime \in (Q-q,Q]$, then $q^\prime -\bar{a}$ is the
unique multiple of $q$ in the interval $(Q-q-\bar{a},Q-\bar{a}]$.
Hence $q^\prime -\bar{a}=q[\frac{Q-\bar{a}}{q}]$, and so we get
\begin{equation}\label{2.1}
q^\prime =q\bigg[\frac{Q-\bar{a}}{q} \bigg]+\bar{a} .
\end{equation}
Making use of $a^{\prime \prime}q-aq^{\prime \prime}=1$ and of
$q^{\prime \prime} \in (Q-q,Q]$, we arrive in a similar way at
\begin{equation}\label{2.2}
q^{\prime \prime}=q\bigg[ \frac{Q+\bar{a}}{q}\bigg]-\bar{a} .
\end{equation}

We consider the partition $I_q^{(1)} \cup I_q^{(2)} \cup I_q^{(3)}
\cup I_q^{(4)}$ of $[0,q]$, where
\begin{equation*}
\begin{split}
& I_q^{(1)}=[\max (2q-Q,0),\min (Q-q,q)] ;\qquad I_q^{(2)} =[\max
(2q-Q,Q-q),q] ;\\
& I_q^{(3)} =[0,\min (2q-Q-1,Q-q)];\qquad \qquad I_q^{(4)}
=[Q-q,2q-Q] .
\end{split}
\end{equation*}

It is clear that $I_q^{(1)}=\emptyset$ unless $q\leq
\frac{2Q}{3}$, $I_q^{(2)}=I_q^{(3)}=\emptyset$ unless $q\geq
\frac{Q}{2}$, and $I_q^{(4)} =\emptyset$ unless $q\geq
\frac{2Q}{3}$. Taking into account \eqref{2.1} and \eqref{2.2} we
see that
\begin{equation*}
q^\prime >q \ \Longleftrightarrow \ \bar{a}\leq Q-q
\end{equation*}
and
\begin{equation*}
q^{\prime \prime} >q \ \Longleftrightarrow \ \bar{a} \geq 2q-Q.
\end{equation*}
Since $\bar{a}$ cannot actually take the values $0$ and $q$, one
has

\medskip

\begin{lem}\label{L2.2}
{\em (i)} $\ \displaystyle \min(q^\prime,q^{\prime \prime})>q \
\Longleftrightarrow \ q\leq \frac{2Q}{3} \ \mbox{and} \ \bar{a}
\in I_q^{(1)}$.

{\em (ii)} $\ \displaystyle q^\prime <q<q^{\prime \prime} \
\Longleftrightarrow \ q\geq \frac{Q}{2} \ \mbox{and} \ \bar{a} \in
I_q^{(2)}$.

{\em (iii)} $\ \displaystyle q^{\prime \prime} <q<q^\prime \
\Longleftrightarrow \ q\geq \frac{Q}{2} \ \mbox{and}\ \bar{a} \in
I_q^{(3)}$.

{\em (iv)} $\ \displaystyle q>\max (q^\prime,q^{\prime \prime}) \
\Longleftrightarrow \ q\geq \frac{2Q}{3} \ \mbox{and} \ \bar{a}
\in I_q^{(4)} .$
\end{lem}

\medskip

For each $\omega$ we put $e_{\omega,Q}=(q,a)$ if the half-line
$\R_+ \omega$ first intersects the scatterer
$(q,a)+V_{1/Q}$ among the $N_Q$ components of
$\mathcal{V}_Q$. We denote by $\omega_{q,a}$ the angle determined by
the trajectories which end near the lattice point $(q,a)$, that is
\begin{equation*}
\omega_{q,a}=\bigg| \bigg\{ \omega \in \Big[ 0,\frac{\pi}{4} \Big)
\, ;\, e_{\omega,Q}=(q,a)\bigg\} \bigg| .
\end{equation*}
\begin{figure}[ht]
\includegraphics*[scale=0.72, bb=0 0 250 135]{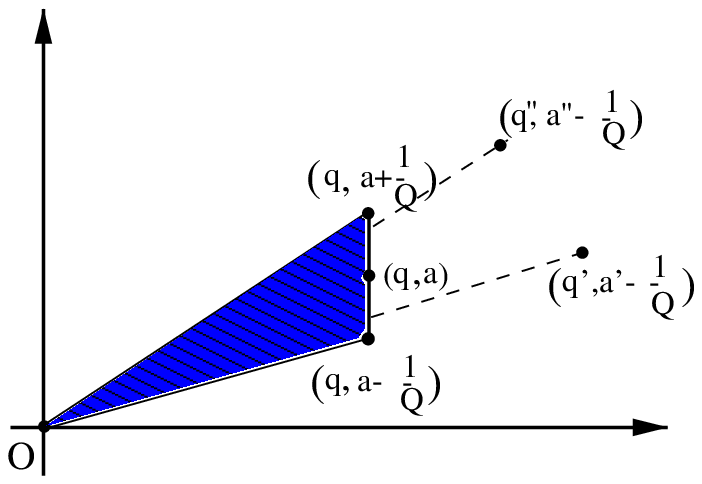}
\includegraphics*[scale=0.72, bb=0 0 230 135]{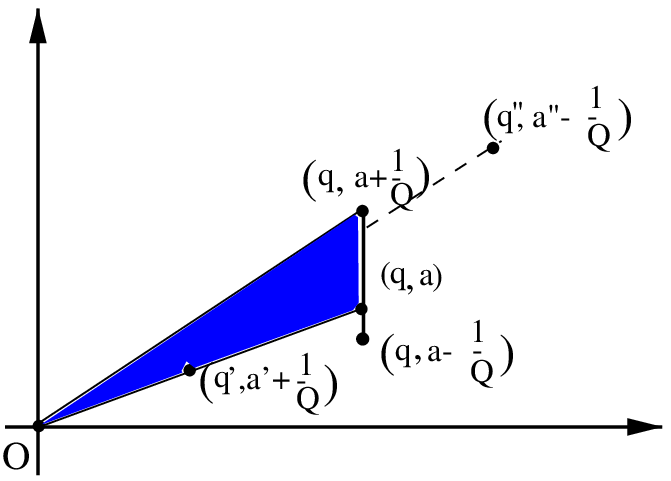}
\caption{The cases $q<\min (q^\prime,q^{\prime \prime})$ and
$q^\prime <q<q^{\prime \prime}$}\label{Figure5}
\end{figure}
\begin{figure}[ht]
\begin{center}
\includegraphics*[scale=0.72, bb=0 0 190 140]{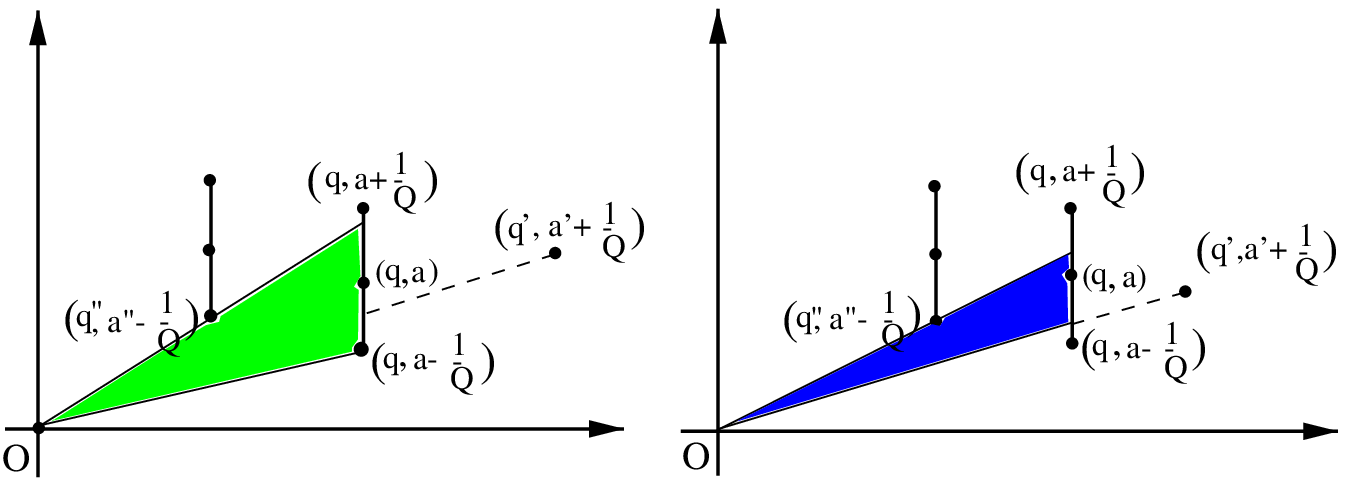}
\includegraphics*[scale=0.72, bb=-40 0 180 140]{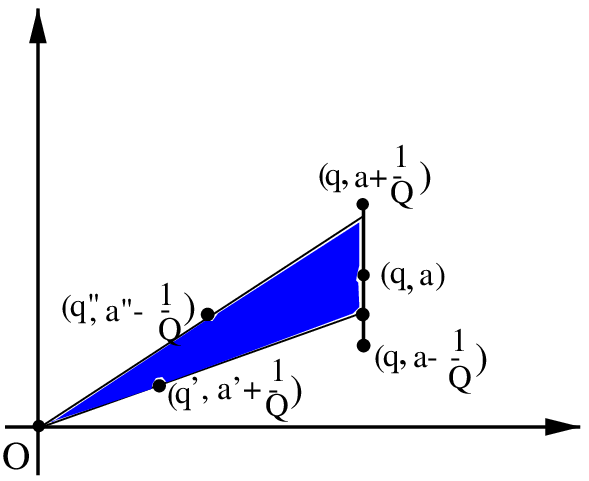}
\end{center}
\caption{The cases $q^{\prime \prime}<q<q^\prime$ and
$q>\max (q^\prime,q^{\prime \prime})$} \label{Figure6}
\end{figure}
\bigskip

Explicit formulas for $\omega_{q,a}$ can be given using Lemma
\ref{L2.2} and
\begin{equation}\label{2.3}
\begin{split}
\arctan (x+h) & -\arctan x=\frac{h}{1+x^2} +O(h^2) \\
& =\frac{h}{1+(x+h)^2} +O(h^2) ,\qquad x\in [0,1],\ h>0 \
\mbox{\rm small} ,
\end{split}
\end{equation}
as follows

\medskip

\begin{lem}\label{L2.3}
{\em (i)} If $\ q\leq \frac{Q}{2}$, then $I_q^{(1)}=[0,q]$,
$I_q^{(2)}=I_q^{(3)}=I_q^{(4)}=\emptyset$, and
\begin{equation*}
\omega_{q,a}=\frac{2}{Qq(1+\gamma^2)}+O\bigg( \frac{1}{Q^2 q}
\bigg) .
\end{equation*}

{\em (ii)} If $\ \frac{Q}{2}<q\leq \frac{2Q}{3}$, then
$I_q^{(4)}=\emptyset$ and
\begin{equation*}
\omega_{q,a}=\begin{cases} \displaystyle \frac{2}{Qq(1+\gamma^2)}
+O\left( \frac{1}{Q^2 q} \right) & \mbox{if $\bar{a}\in I_q^{(1)}
=[2q-Q,Q-q]$;} \\ & \\ \displaystyle \frac{Q-q+\bar{a}}{Qq\bar{a}
(1+\gamma^2)}+O\left( \frac{1}{Q^2 q^\prime} \right) & \mbox{if
$\bar{a}\in I_q^{(2)} =[Q-q,q]$;}
\\ & \\ \displaystyle \frac{Q-\bar{a}}{Qq(q-\bar{a})(1+\gamma^2)}
+O\left( \frac{1}{Q^2 q^{\prime \prime}} \right) & \mbox{if
$\bar{a} \in I_q^{(3)}=[0,2q-Q]$.}
\end{cases}
\end{equation*}

{\em (iii)} If $\ q>\frac{2Q}{3}$, then $I_q^{(1)}=\emptyset$ and
\begin{equation*}
\omega_{q,a} =\begin{cases} \displaystyle
\frac{Q-q+\bar{a}}{Qq\bar{a}(1+\gamma^2)} +O\left( \frac{1}{Q^2
q^\prime} \right) & \mbox{if $\bar{a}\in I_q^{(2)}=[2q-Q,q]$;}
\\ & \\ \displaystyle \frac{Q-\bar{a}}{Qq(q-\bar{a})(1+\gamma^2)}
+O\left( \frac{1}{Q^2 q^{\prime \prime}} \right) & \mbox{if
$\bar{a}\in I_q^{(3)}=[0,Q-q]$;} \\ & \\ \displaystyle
\frac{Q-q}{Q\bar{a} (q-\bar{a})(1+\gamma^2)}+O\bigg( \frac{1}{Q^2
q^\prime}+\frac{1}{Q^2 q^{\prime \prime}} \bigg) & \mbox{if
$\bar{a}\in I_q^{(4)}=[Q-q,2q-Q]$.}
\end{cases}
\end{equation*}
\end{lem}

{\sl Proof.} If $q>\frac{Q}{2}$ then one has, according to Lemma
\ref{L2.2}, that $q<\min (q^\prime,q^{\prime \prime})$. In this
case we infer from \eqref{2.3} that
\begin{equation*}
\omega_{q,a} =\arctan \frac{a+\frac{1}{Q}}{q}-\arctan
\frac{a-\frac{1}{Q}}{q} =\frac{2}{Qq \Big(
1+\frac{(a-\frac{1}{Q})^2}{q^2} \Big)} +O\bigg( \frac{1}{Q^2 q^2}
\bigg) .
\end{equation*}
This implies, in connection with
\begin{equation*}
\frac{1}{1+\frac{(a-\frac{1}{Q})^2}{q^2}}-\frac{1}{1+\frac{a^2}{q^2}}
\ll \frac{a^2}{q^2} -\frac{(a-\frac{1}{Q})^2}{q^2} \ll \frac{2
a}{Q q} \ll \frac{1}{Q} \, ,
\end{equation*}
the desired formula for $\omega_{q,a}$.

In (ii) and (iii) we split the discussion in the following three
cases:

1) $\ \bar{a} \in I_q^{(2)}$, that is $q^\prime <q<q^{\prime
\prime}$.

In this case
$\frac{a-\frac{1}{Q}}{q}<\frac{a^\prime+\frac{1}{Q}}{q^\prime}
<\frac{a}{q} <\frac{a^{\prime \prime}-\frac{1}{Q}}{q^{\prime
\prime}} <\frac{a+\frac{1}{Q}}{q}$ and we may write
\begin{equation*}
\begin{split}
\omega_{q,a} & =\arctan \frac{a+\frac{1}{Q}}{q}-\arctan
\frac{a^\prime +\frac{1}{Q}}{q^\prime}
=\frac{1-\frac{q-q^\prime}{Q}}{qq^\prime \Big( 1+\frac{(a^\prime
+\frac{1}{Q})^2}{q^{\prime 2}} \Big)} +O\Bigg( \bigg(
\frac{(1-\frac{q-q^\prime}{Q})}{qq^\prime}\bigg)^2 \Bigg) \\
& =\frac{Q-q+q^\prime}{Qqq^\prime (1+\frac{a^{\prime 2}}{q^{\prime
2}} )} +O\bigg( \frac{1}{Qq^2
q^\prime}+\frac{(Q-q+q^\prime)^2}{Q^2 q^2 q^{\prime 2}} \bigg) \\
& =\frac{Q-q+q^\prime}{Qqq^\prime (1+\gamma^2)} +O \bigg(
\frac{Q-q+q^\prime}{qq^\prime} \Big(
\frac{a}{q}-\frac{a^\prime}{q^\prime} \Big) +\frac{1}{Qq^2
q^\prime}+\frac{(Q-q+q^\prime)^2}{Q^2 q^2 q^{\prime 2}} \bigg) .
\end{split}
\end{equation*}
The desired formula for $\omega_{q,a}$ follows in this case from
the formula above, the fact that $q^\prime =\bar{a}$ as a result
of $[\frac{Q-\bar{a}}{q} ]=0$, from $q\geq \frac{Q}{2}$, and from
\begin{equation*}
\begin{split}
& \frac{Q-q+q^\prime}{qq^\prime} \bigg(
\frac{a}{q}-\frac{a^\prime}{q^\prime} \bigg)
=\frac{Q-q+q^\prime}{q^2 q^{\prime 2}} <\frac{2q^\prime}{q^2
q^{\prime 2}} =\frac{2}{q^2 q^\prime} \ll \frac{1}{Q^2 q^\prime}
\, ,\\
& \frac{1}{Qq^2 q^\prime} \ll \frac{1}{Q^3 q^\prime} \leq
\frac{1}{Q^2 q^\prime} \, ,\\
& \frac{(Q-q+q^\prime)^2}{Q^2 q^2 q^{\prime 2}} \ll \frac{1}{q^2
q^{\prime 2}} \ll \frac{1}{Q^2 q^{\prime 2}} \leq \frac{1}{Q^2
q^\prime} \, .
\end{split}
\end{equation*}

\smallskip

2) $\ \bar{a}\in I_q^{(3)}$, that is $q^{\prime \prime}
<q<q^\prime$.

In this case $q^{\prime \prime}=q-\bar{a}$, $q\geq \frac{Q}{2}$,
and $\frac{a-\frac{1}{Q}}{q}<\frac{a^\prime+\frac{1}{Q}}{q^\prime}
<\frac{a}{q}<\frac{a^{\prime \prime}-\frac{1}{Q}}{q^{\prime
\prime}} <\frac{a+\frac{1}{Q}}{q}$. Proceeding as in the case
$\bar{a}\in I_q^{(2)}$ we find
\begin{equation*}
\begin{split}
\omega_{q,a} & =\arctan \frac{a^{\prime
\prime}-\frac{1}{Q}}{q^{\prime \prime}} -\arctan
\frac{a-\frac{1}{Q}}{q} =\frac{1-\frac{q-q^{\prime
\prime}}{Q}}{qq^{\prime \prime} \Big(
1+\frac{(a-\frac{1}{Q})^2}{q^2} \Big)} +O\left( \frac{1}{Q^2
q^{\prime \prime 2}} \right) \\
& =\frac{Q-\bar{a}}{Qq(q-\bar{a})(1+\gamma^2)} +O\left(
\frac{1}{Q^2 q^{\prime \prime}} \right),
\end{split}
\end{equation*}
as required.

\smallskip

3) $\ \bar{a}\in I_q^{(4)}$, that is $q>\max (q^\prime ,q^{\prime
\prime})$.

In this case
$\frac{a-\frac{1}{Q}}{q}<\frac{a^\prime+\frac{1}{Q}}{q^\prime}
<\frac{a}{q}<\frac{a^{\prime \prime}-\frac{1}{Q}}{q^{\prime
\prime}} <\frac{a+\frac{1}{Q}}{q}$, and we write
\begin{equation*}
\omega_{q,a}=\arctan \frac{a^{\prime
\prime}-\frac{1}{Q}}{q^{\prime \prime}}-\arctan \frac{a^\prime
+\frac{1}{Q}}{q^\prime} =\omega_{q,a}^{(1)}+\omega_{q,a}^{(2)} ,
\end{equation*}
with
\begin{equation*}
\begin{split}
& \omega_{q,a}^{(1)} =\arctan \frac{a^{\prime
\prime}-\frac{1}{Q}}{q^{\prime \prime}} -\arctan \frac{a}{q}=
\frac{Q-q}{Qqq^{\prime \prime} (1+\gamma^2)}+O\left( \frac{1}{Q^2
q^{\prime \prime 2}} \right) ,\\ & \omega_{q,a}^{(2)} =\arctan
\frac{a}{q}-\arctan \frac{a^\prime+\frac{1}{Q}}{q^\prime}
=\frac{Q-q}{Qqq^\prime (1+\gamma^2)} +O\left( \frac{1}{Q^2
q^{\prime 2}}\right) .
\end{split}
\end{equation*}
Since in this case $q^{\prime \prime}=q-\bar{a}$ and $q^\prime
=\bar{a}$ (thus $q^\prime +q^{\prime \prime} =q$) we arrive at
\begin{equation*}
\omega_{q,a} =\frac{Q-q}{Q\bar{a}(q-\bar{a})(1+\gamma^2)}+O\bigg(
\frac{1}{Q^2 q^{\prime 2}}+\frac{1}{Q^2 q^{\prime \prime 2}}
\bigg) ,
\end{equation*}
as required. \qed

\bigskip

\section{The case of vertical scatterers}
In this section we estimate $\tilde{H}_{\eps,I,Q}(t)$ as $\eps
\rightarrow 0^+$ in the hypothesis that $\vert I\vert \rightarrow
0$ and $Q\rightarrow \infty$ in a controlled way. More precisely
we prove

\medskip

\begin{prop}\label{P3.1}
Let $\theta,\theta_1 \in (0,1)$ and suppose that $I=[\tan
\omega_0,\tan \omega_1]$ is a subinterval of $[0,1]$ of size
$\vert I\vert \asymp \eps^\theta$, and $Q$ is an integer such that
$Q=\frac{\cos \omega_0}{\eps}+O(\eps^{-\theta})$. Then one has
\begin{equation*}
\tilde{H}_{\eps,I,Q}(t)=c_I H(t)+O\big(
\EE_{\theta,\theta_1,\delta}(\eps)\big) \qquad \mbox{as $\ \eps
\rightarrow 0^+$},
\end{equation*}
uniformly for $t$ in a compact subset of $(0,\infty )\setminus
\{ 1,2\}$, where $H$ has been explicitly defined in Theorem
\ref{T1.1}, and we denote
\begin{equation*}
c_I=\int\limits_I \frac{dt}{1+t^2} =\omega_1-\omega_0,\qquad
\EE_{\theta,\theta_1,\delta}(\eps) =
\eps^{\min(2\theta,\frac{1}{2}-2\theta_1-\delta,\theta+\theta_1-\delta)}.
\end{equation*}
\end{prop}

\medskip

It is convenient to write
\begin{equation}\label{0.1}
\tilde{H}_{\eps,I,Q}(t)=\bigg| \bigg\{ \omega \, ;\, \tan \omega
\in I,\ l_{\frac{1}{Q}} (\omega) \cos \omega
>\frac{t\cos \omega}{\eps} \bigg\} \bigg| .
\end{equation}

We also define for each interval $I\subseteq [0,1]$ and each
integer $Q\geq 1$ the quantity
\begin{equation*}
\tilde{H}_{I,Q}(t)=\bigg| \bigg\{ \omega \, ;\, \tan \omega \in
I,\ l_{\frac{1}{Q}} (\omega)\cos \omega >tQ \bigg\} \bigg| .
\end{equation*}

Let $Q^-$ and $Q^+$ be two integers such that
\begin{equation*}
Q^-\leq \frac{\cos \omega_1}{\eps} \leq \frac{\cos \omega_0}{\eps}
\leq Q^+ .
\end{equation*}
This also gives for any $Q$ as in the statement of Proposition
\ref{P3.1}
\begin{equation}\label{0.2}
\frac{Q^\pm}{Q}-1=O(\eps^\theta ).
\end{equation}

Since $\eps \mapsto \tilde{H}_{\eps,I,Q}(t)$ is monotonically
increasing, one has
\begin{equation}\label{0.3}
\tilde{H}_{I,Q} \bigg( \frac{tQ^+}{Q} \bigg) \leq
\tilde{H}_{\eps,I,Q}(t) \leq \tilde{H}_{\eps,I,Q} \bigg(
\frac{tQ^-}{Q} \bigg) .
\end{equation}

We now estimate $\tilde{H}_{I,Q}(t)$ in

\medskip

\begin{prop}\label{P3.2}
Let $\theta,\theta_1 \in (0,1)$ and suppose that $I=[\tan
\omega_0,\tan \omega_1]$ is a subinterval of $[0,1]$ of size
$\vert I\vert \asymp Q^{-\theta}$. Then one has
\begin{equation*}
\tilde{H}_{I,Q}(t)=c_I H(t)+O\big(
E_{\theta,\theta_1,\delta}(Q)\big) \qquad \mbox{as $\ Q\rightarrow
\infty$},
\end{equation*}
where we denote
\begin{equation*}
E_{\theta,\theta_1,\delta} (Q)=Q^{\max
(2\theta_1-\frac{1}{2}+\delta,-\theta-\theta_1+\delta)}.
\end{equation*}
\end{prop}

\medskip

The key technical tools we shall employ to estimate
$\tilde{H}_{I,Q}(t)$ are the following two lemmas from \cite{BGZ}
and \cite{BCZ}.

\medskip

\begin{lem}\label{L3.3}
{\em (\cite[Lemma 2.2]{BGZ})} Suppose that $q\geq 1$ is an
integer, $\II$ and $\JJ$ are intervals of length lesser than $q$,
and that $f$ is a $C^1$ function on $\II \times \JJ$. Then one has
\begin{equation*}
\begin{split}
\sum\limits_{\substack{a\in \II,\, b\in \JJ \\ ab=1\hspace{-6pt}
\pmod{q}}} \hspace{-8pt} f(a,b) & =\frac{\varphi (q)}{q^2}
\iint\limits_{\II \times \JJ} f(x,y)\, dx\, dy \\
& +O_\delta \left( T^2 q^{\frac{1}{2}+\delta} \| f\|_\infty +
Tq^{\frac{3}{2}+\delta} \| Df\|_\infty +\frac{\vert \II \vert \,
\vert \JJ \vert \, \| Df\|_\infty}{T} \right)
\end{split}
\end{equation*}
for any integer $T\geq 1$ and any $\delta >0$, where $\| \cdot\|$
denotes the $L^\infty$ norm on $\II \times \JJ$, $Df=\vert
\frac{\partial f}{\partial x} \vert+\vert \frac{\partial
f}{\partial y} \vert$, and $\varphi$ is Euler's totient function.
\end{lem}

\medskip

\begin{lem}\label{L3.4}
{\em (\cite[Lemma 2.3]{BCZ})} Suppose that $0<a<b$ and that $f$ is
a $C^1$ function on $[a,b]$. Then one has
\begin{equation*}
\sum\limits_{a<k\leq b} \frac{\varphi(k)}{k}\,
f(k)=\frac{1}{\zeta(2)} \int\limits_a^b f(x)\, dx+O\Bigg( \ln b
\bigg( \| f\|_\infty +\int\limits_a^b \vert f^\prime \vert \bigg)
\Bigg) .
\end{equation*}
\end{lem}

\medskip

We shall need the next two corollaries of these two lemmas.

\medskip

\begin{cor}\label{C3.5}
{\em (i)} For each $0<t_1<t_2 \leq 1$ one has
\begin{equation*}
 \sum\limits_{t_1 Q<q\leq t_2 Q} \, \sum\limits_{\substack{a\in qI,\,
\bar{a} \in [0,q] \\ a\bar{a} =1\hspace{-6pt} \pmod{q}}}
\hspace{-2pt} \frac{1}{Qq(1+\gamma^2)} =\frac{c_I
(t_2-t_1)}{\zeta(2)}+O \big( E_{\theta,\theta_1,\delta}(Q)\big).
\end{equation*}

{\em (ii)} For each $\frac{1}{2}<t_1<t_2 \leq 1$ one has
\begin{equation*}
\begin{split}
& \sum\limits_{t_1 Q<q\leq t_2 Q} \ \sum\limits_{\substack{a\in qI
,\, \bar{a} \in [0,Q-q] \\ a\bar{a} =1\hspace{-6pt} \pmod{q}}}
\frac{1}{Qq(1+\gamma^2)} =\frac{c_I}{\zeta(2)}
\int\limits_{t_1}^{t_2} \frac{1-x}{x}\, dx
+O \big( E_{\theta,\theta_1,\delta}(Q)\big);\\
& \sum\limits_{t_1 Q<q\leq t_2 Q} \ \sum\limits_{\substack{a\in
qI,\, \bar{a} \in [2q-Q,q] \\ a\bar{a} =1\hspace{-6pt} \pmod{q}}}
\frac{1}{Qq(1+\gamma^2)} =\frac{c_I}{\zeta(2)}
\int\limits_{t_1}^{t_2} \frac{1-x}{x}\, dx+O \big(
E_{\theta,\theta_1,\delta}(Q)\big) .
\end{split}
\end{equation*}

{\em (iii)} For each $t\in (\frac{1}{2},\frac{2}{3})$ one has
\begin{equation*}
\sum\limits_{tQ<q\leq \frac{2Q}{3}} \sum\limits_{\substack{a\in qI \\
\bar{a} \in [2q-Q,Q-q] \\ a\bar{a} =1\hspace{-6pt} \pmod{q}}}
\hspace{-6pt} \frac{1}{Qq(1+\gamma^2)} =c_I
\int\limits_t^{2/3} \frac{2-3x}{x}\, dx+O \big(
E_{\theta,\theta_1,\delta}(Q)\big) .
\end{equation*}
\end{cor}

\medskip

\begin{cor}\label{C3.6}
For each $\frac{1}{2}<t_1<t_2 \leq 1$, both sums
\begin{equation*}
\begin{split}
& \sum\limits_{t_1 Q<q\leq t_2 Q} \ \sum\limits_{\substack{a\in qI,\,
\bar{a}\in (Q-q,q] \\ a\bar{a}=1 \hspace{-6pt} \pmod{q}}}
\frac{Q-q}{Qq\bar{a}(1+\gamma^2)} \qquad \mbox{and}
\\ &
\sum\limits_{t_1 Q<q\leq t_2 Q} \ \sum\limits_{\substack{a\in qI ,\,
\bar{a}\in [0,2q-Q) \\ a\bar{a}=1 \hspace{-6pt} \pmod{q}}}
\frac{Q-q}{Qq(q-\bar{a})(1+\gamma^2)}
\end{split}
\end{equation*}
are of the form
\begin{equation*}
\frac{c_I}{\zeta(2)} \int\limits_{t_1}^{t_2} \frac{1-x}{x} \, \ln
\frac{x}{1-x} \, dx+O \big( E_{\theta,\theta_1,\delta}(Q)\big).
\end{equation*}
\end{cor}

\medskip

We only give the proof of Corollary \ref{C3.6}. The proof of
Corollary \ref{C3.5} is easier and we leave it as an exercise to
the reader.

\medskip

{\sl Proof of Corollary \ref{C3.6}.} We only prove the first
equality. The second one is proved in a similar way by changing
$\bar{a}$ to $q-\bar{a}$. To estimate the inner sum we consider
for each $q\in (t_1 Q,t_2 Q]\subset (\frac{Q}{2},Q]$ the function
\begin{equation*}
f(a,\bar{a})=f_q(a,\bar{a})=\frac{Q-q}{Qq\bar{a}(1+\gamma^2)}
=\frac{q(Q-q)}{Q(q^2+a^2)\bar{a}} \, ,\quad a\in qI,\ \bar{a}\in
[Q-q+1,q].
\end{equation*}
Since $Q-q+1\leq \bar{a}$, one has
\begin{equation*}
\| f\|_\infty \leq \frac{1}{Qq}\, ,\quad \left\| \frac{\partial
f}{\partial a} \right\|_\infty \ll \frac{1}{Qq^2} \quad \mbox{\rm
and} \quad \left\| \frac{\partial f}{\partial \bar{a}}
\right\|_\infty \ll \frac{1}{Qq(Q-q+1)} \, .
\end{equation*}

Applying Lemma \ref{L3.3} to $f$ with $T=[Q^{\theta_1}]$,
$\II=[Q-q+1,q]$ and $\JJ=qI$, we gather
\begin{equation}\label{I}
\sum\limits_{\substack{a\in qI,\, \bar{a}\in (Q-q,q] \\ a\bar{a}=1
\hspace{-6pt} \pmod{q}}}  \frac{Q-q}{Qq\bar{a}(1+\gamma^2)}
=\frac{\varphi(q)}{q} \cdot W_Q (q)+O\big( F_Q(q)\big) ,
\end{equation}
where
\begin{equation*}
W_Q(q) =\frac{Q-q}{Qq^2} \int\limits_{qI}
\frac{da}{1+\frac{a^2}{q^2}} \int\limits_{Q-q+1}^q
\frac{d\bar{a}}{\bar{a}}=\frac{c_I (Q-q)}{Qq} \, \ln
\frac{q}{Q-q+1} \, ,
\end{equation*}
and the error terms
\begin{equation*}
\begin{split}
F_Q(q) & =Q^{2\theta_1-1} q^{-\frac{1}{2}+\delta} +Q^{\theta_1 -1}
\frac{q^{\frac{1}{2}+\delta}}{Q-q+1} +\frac{qQ^{-\theta}
(2q-Q)\frac{1}{Qq(Q-q+1)}}{Q^{\theta_1}}
\\ & \leq Q^{2\theta_1-\frac{3}{2}+\delta}
+\frac{Q^{\theta_1-\frac{1}{2}+\delta}}{Q-q+1}
+\frac{Q^{-\theta-\theta_1}}{Q-q+1}
\end{split}
\end{equation*}
sum up over $q$ to (we may take $\delta <\theta_1$)
\begin{equation}\label{II}
\sum\limits_{\frac{Q}{2}<q\leq Q} F_Q(q)  \ll
Q^{2\theta_1-\frac{1}{2}+\delta} +Q^{\theta_1-\frac{1}{2}+\delta}
\ln Q+Q^{-\theta-\theta_1} \ln Q \ll
F_{\theta,\theta_1,\delta}(Q).
\end{equation}

Using the inequality $\ln (1+x)\leq x$, we see immediately that
\begin{equation}\label{III}
\begin{split}
\int\limits_{t_1 Q}^{t_2 Q} W_Q(q)\, dq & =c_I \int\limits_{t_1
Q}^{t_2 Q} \frac{Q-q}{Qq} \, \ln \frac{q}{Q-q} \, dq+O\left(
\frac{1}{Q} \right) \\
& =c_I \int\limits_{t_1}^{t_2} \frac{1-x}{x} \, \ln
\frac{x}{1-x}\, dx+O(Q^{-1} ) .
\end{split}
\end{equation}

Finally, the inequalities $\| W_Q \|_{\infty,[t_1 Q,t_2
Q]} \ll Q^{\delta -1}$ and $\| W_Q^\prime \|_{\infty ,[t_1 Q,t_2
Q]} \ll q^{\delta-2}$ show that Lemma \ref{L3.4} applied to the
function $W_Q$ yields
\begin{equation}\label{IV}
\sum\limits_{t_1 Q<q\leq t_2 Q} \frac{\varphi(q)}{q}\, W_Q(q)
=\frac{1}{\zeta(2)} \int\limits_{t_1 Q}^{t_2 Q} W_Q(q)\,
dq+O(Q^{2\delta-1}) .
\end{equation}

Finally we put together \eqref{I}, \eqref{II}, \eqref{III} and
\eqref{IV} to get the desired estimate in the first formula of the
statement. \qed

\medskip

{\sl Proof of Proposition \ref{P3.2}.} We use the explicit
formulas for $\omega_{q,a}$ found in Section 2. The crude estimate
\begin{equation*}
\begin{split}
& \max \bigg( \sum\limits_{\frac{a}{q}\in \FI} \frac{1}{Q^2 q} \,
, \sum\limits_{\frac{a}{q}\in \FI} \frac{1}{Q^2 q^\prime}\, ,
\sum\limits_{\frac{a}{q}\in \FI} \frac{1}{Q^2 q^{\prime \prime}}
\bigg) \leq \frac{1}{Q^2} \sum\limits_{q=1}^Q \frac{\# (qI\cap
\Z)+1}{q} \\
& \qquad \ll \frac{1}{Q^2} \sum\limits_{q=1}^Q \frac{q\vert
I\vert+1}{q} \ll \ln Q+\frac{\vert I\vert}{Q} \ll Q^{-\theta-1}
\end{split}
\end{equation*}
and Lemma \ref{L2.3} show that
\begin{equation*}
\tilde{H}_{I,Q}(t)=S_{I,Q}(t)+O(Q^{-\theta-1}) ,
\end{equation*}
where
\begin{equation*}
S_{I,Q}(t)=\sum\limits_{\substack{\frac{a}{q}\in \FI \\ q>tQ}}
\tilde{\omega}_{q,a} ,
\end{equation*}
with
\begin{equation}\label{3.1}
\tilde{\omega}_{q,a} =\begin{cases} \displaystyle
\frac{2}{Qq(1+\gamma^2)} & \mbox{if $\bar{a}\in I_q^{(1)}$;} \\
&
\\ \displaystyle \frac{Q-q+\bar{a}}{Qq\bar{a}(1+\gamma^2)} &
\mbox{if $\bar{a}\in I_q^{(2)}$;}
\\ & \\ \displaystyle \frac{Q-\bar{a}}{Qq(q-\bar{a})(1+\gamma^2)}
& \mbox{if $\bar{a}\in I_q^{(3)}$;}
\\ & \\ \displaystyle \frac{Q-q}{Q\bar{a}
(q-\bar{a})(1+\gamma^2)} & \mbox{if $\bar{a}\in I_q^{(4)}$.}
\end{cases}
\end{equation}
It is clear that
\begin{equation}\label{3.2}
\sum\limits_{\frac{a}{q} \in \FI} \tilde{\omega}_{q,a}
=\int\limits_{\omega_0}^{\omega_1} d\omega +O\left(
\frac{1}{Q}\right) =c_I+O(Q^{-1}) .
\end{equation}

Suppose now that $0<t<\frac{1}{2}$. In this case we may write in
view of Lemma \ref{L2.3}, \eqref{3.1}, \eqref{3.2} and Corollary
\ref{C3.5} (i)
\begin{equation*}
\begin{split}
S_{I,Q}(t) & =\sum\limits_{\frac{a}{q} \in \FI}
\tilde{\omega}_{q,a} -\sum\limits_{\substack{\frac{a}{q}\in \FI \\
q\leq tQ}} \tilde{\omega}_{q,a} \\
& =c_I -\sum\limits_{0<q\leq tQ} \sum\limits_{\substack{a\in qI,\,
\bar{a} \in [0,q] \\ a\bar{a} =1\hspace{-6pt} \pmod{q}}}
\hspace{-2pt} \frac{2}{Qq(1+\gamma^2)} +O(Q^{-1}) \\
& =c_I-\frac{2c_I t}{\zeta(2)} +O \big( E_{\theta,\theta_1,\delta}
(Q)\big) .
\end{split}
\end{equation*}

If $t\geq \frac{2}{3}$, then Lemma \ref{L2.3} and \eqref{3.1}
yield
\begin{equation*}
\begin{split}
S_{I,Q}(t) & =\sum\limits_{tQ<q\leq Q} \Bigg(
\sum\limits_{\substack{a\in qI ,\,\bar{a} \in I_q^{(2)} \\
a\bar{a}=1\hspace{-6pt}\pmod{q}}} \hspace{-8pt}
\tilde{\omega}_{q,a}+\sum\limits_{\substack{a\in qI ,\, \bar{a}
\in I_q^{(3)} \\ a\bar{a}=1\hspace{-6pt}\pmod{q}}} \hspace{-8pt}
\tilde{\omega}_{q,a}+\sum\limits_{\substack{a\in qI ,\,\bar{a} \in
I_q^{(4)}\\ a\bar{a}=1\hspace{-6pt}\pmod{q}}} \hspace{-8pt}
\tilde{\omega}_{q,a} \Bigg) \\
& =\sum\limits_{tQ<q\leq Q} \ \sum\limits_{\substack{a\in
qI,\, \bar{a}\in [2q-Q,q] \\ a\bar{a}=1\hspace{-6pt}\pmod{q}}}
\frac{Q-q+\bar{a}}{Qq\bar{a}(1+\gamma^2)} \\
& \qquad +\sum\limits_{tQ<q\leq Q} \
\sum\limits_{\substack{a\in qI ,\,\bar{a}\in [0,Q-q] \\
a\bar{a}=1\hspace{-6pt} \pmod{q}}} \hspace{-2pt}
\frac{Q-\bar{a}}{Qq(q-\bar{a})(1+\gamma^2)} \\
& \qquad \qquad +\sum\limits_{tQ<q\leq Q}\
\sum\limits_{\substack{a\in qI \\ \bar{a}\in
(Q-q,2q-Q)
\\ a\bar{a}=1\hspace{-6pt} \pmod{q}}} \hspace{-2pt}
\frac{Q-q}{Q\bar{a}(q-\bar{a})(1+\gamma ^2)} \, ,
\end{split}
\end{equation*}
which can also be written as
\begin{equation*}
\begin{split}
& \sum\limits_{tQ<q\leq Q} \ \sum\limits_{\substack{a\in qI ,\,
\bar{a}\in (Q-q,q] \\ a\bar{a} =1\hspace{-6pt} \pmod{q}}}
\frac{Q-q}{Qq\bar{a}(1+\gamma^2)} \\
& \qquad +\sum\limits_{tQ<q\leq Q} \
\sum\limits_{\substack{a\in qI ,\, \bar{a} \in [0,2q-Q) \\
a\bar{a}=1\hspace{-6pt} \pmod{q}}} \hspace{-2pt}
\frac{Q-q}{Qq(q-\bar{a})(1+\gamma^2)} \\
& \qquad \qquad +\sum\limits_{tQ<q\leq Q} \
\sum\limits_{\substack{a\in qI \\
\bar{a} \in [0,Q-q] \cup [2q-Q,q] \\
a\bar{a}=1\hspace{-6pt}\pmod{q}}} \hspace{-8pt}
\frac{1}{Qq(1+\gamma^2)} \, .
\end{split}
\end{equation*}

Applying Corollaries \ref{C3.6} and \ref{C3.5} (ii) we find that
\begin{equation}\label{3.3}
\begin{split}
S_{I,Q}(t) & =\frac{2c_I}{\zeta(2)} \int\limits_t^1 \frac{1-x}{x}
\left( 1+\ln \frac{x}{1-x}\right) dx+ O\big(
E_{\theta,\theta_1,\delta}(Q)\big) \\ & =\frac{2c_I}{\zeta(2)}
\int\limits_t^1 \psi(x)\, dx+O\big(
E_{\theta,\theta_1,\delta}(Q)\big) .
\end{split}
\end{equation}

Finally, in the case where $\frac{1}{2}<t<\frac{2}{3}$, Lemma
\ref{L2.3} yields
\begin{equation}\label{3.4}
S_{I,Q}(t)=\sum\limits_{tQ<q\leq Q} \sum\limits_{\substack{a\in qI \\
\gcd (a,q)=1}} \tilde{\omega}_{q,a} =S_{I,Q}\left(
\frac{2}{3}\right)+T_{I,Q}(t),
\end{equation}
where
\begin{equation*}
T_{I,Q}(t)=\sum\limits_{tQ<q\leq \frac{2Q}{3}} \Bigg(
\sum\limits_{\substack{a\in qI ,\,\bar{a} \in I_q^{(1)} \\
a\bar{a}=1\hspace{-6pt}\pmod{q}}} \hspace{-8pt}
\tilde{\omega}_{q,a}+\sum\limits_{\substack{a\in qI ,\, \bar{a}
\in I_q^{(2)} \\ a\bar{a}=1\hspace{-6pt}\pmod{q}}} \hspace{-8pt}
\tilde{\omega}_{q,a}+\sum\limits_{\substack{a\in qI ,\,\bar{a} \in
I_q^{(3)}\\ a\bar{a}=1\hspace{-6pt}\pmod{q}}} \hspace{-8pt}
\tilde{\omega}_{q,a} \Bigg).
\end{equation*}
Using also \eqref{3.1} the sum $T_{I,Q}(t)$ can be successively
written as
\begin{equation*}
\begin{split}
& \mbox{\small $\displaystyle
\sum\limits_{tQ<q\leq \frac{2Q}{3}} \Bigg(
\sum\limits_{\substack{a\in qI \\ \bar{a}\in [2q-Q,Q-q] \\
a\bar{a}=1\hspace{-6pt} \pmod{q}}} \hspace{-8pt}
\frac{2}{Qq(1+\gamma^2)}+\hspace{-8pt}\sum\limits_{\substack{a\in qI \\
\bar{a} \in (Q-q,q] \\ a\bar{a}=1\hspace{-6pt}\pmod{q}}}
\hspace{-8pt} \frac{Q-q+\bar{a}}{Qq\bar{a}(1+\gamma^2)}
+\hspace{-8pt}
\sum\limits_{\substack{a\in qI \\
\bar{a}\in [0,2q-Q) \\
a\bar{a}=1\hspace{-6pt}\pmod{q}}}\hspace{-8pt}
\frac{Q-\bar{a}}{Qq(q-\bar{a})(1+\gamma^2)} \Bigg)$} \\
& \quad \mbox{\small $\displaystyle
 =\sum\limits_{tQ<q\leq \frac{2Q}{3}} \sum\limits_{\substack{a\in
qI \\ \bar{a}\in (Q-q,q] \\ a\bar{a}=1\hspace{-6pt}\pmod{q}}}
\hspace{-8pt} \frac{Q-q}{Qq\bar{a}(1+\gamma^2)}
+\sum\limits_{tQ<q\leq \frac{2Q}{3}} \sum\limits_{\substack{a\in
qI \\ \bar{a} \in [0,2q-Q) \\ a\bar{a}=1\hspace{-6pt}\pmod{q}}}
\hspace{-8pt} \frac{Q-q}{Qq(q-\bar{a})(1+\gamma^2)}$} \\
& \qquad \mbox{\small $\displaystyle +\sum\limits_{tQ<q\leq \frac{2Q}{3}}
\sum\limits_{\substack{a\in qI \\ \bar{a}\in [0,q] \\
a\bar{a}=1\hspace{-6pt} \pmod{q}}} \hspace{-8pt}
\frac{1}{Qq(1+\gamma^2)} +\sum\limits_{tQ<q\leq \frac{2Q}{3}}
\sum\limits_{\substack{a\in qI \\ \bar{a} \in [2q-Q,Q-q] \\
a\bar{a}=1\hspace{-6pt} \pmod{q}}} \hspace{-8pt}
\frac{1}{Qq(1+\gamma^2)} \, .$}
\end{split}
\end{equation*}
Applying Corollaries \ref{C3.5} and \ref{C3.6} we find that
\begin{equation*}
\begin{split}
T_{I,Q}(t) & =\frac{2c_I}{\zeta(2)} \int\limits_t^{2/3}
\frac{1-x}{x} \, \ln \frac{x}{1-x}\, dx+\frac{c_I}{\zeta(2)}
\int\limits_t^{2/3} \left( 1+\frac{2-3x}{x}\right) dx+
O\big( E_{\theta,\theta_1,\delta}(Q)\big)
\\ & =\frac{2c_I}{\zeta(2)} \int\limits_t^{2/3} \psi(x)\,
dx+O\big( E_{\theta,\theta_1,\delta}(Q)\big) .
\end{split}
\end{equation*}
The desired estimate now follows from \eqref{3.3}, \eqref{3.4} and
the equality above. \qed

\medskip

{\sl Proof of Proposition \ref{P3.1}.} By Proposition \ref{P3.2},
the mean value theorem and \eqref{0.2} we infer that
\begin{equation*}
\begin{split}
\tilde{H}_{I,Q} \bigg( \frac{tQ^\pm}{Q} \bigg) & =c_I H\bigg(
\frac{tQ^\pm}{Q}\bigg) +O\big( E_{\theta,\theta_1,\delta}(Q)\big)
\\ & =c_I \big( H(t)+O(\eps^\theta)\big) +O\big(
E_{\theta,\theta_1,\delta} (\theta)\big) \\
& =c_I H(t)+O\big( Q^{-\min
(2\theta,\frac{1}{2}-2\theta_1-\delta,\theta+\theta_1-\delta)}\big)
.
\end{split}
\end{equation*}

The statement in Proposition \ref{P3.1} now follows from this
inequality and \eqref{0.3}. \qed

\medskip

\begin{prop}\label{P3.7}
Let $\theta,\theta_1 \in (0,1)$ and suppose that $I=[\tan
\omega_0,\tan \omega_1]$ is a subinterval of $[0,1]$ of size
$\vert I\vert \asymp \eps^\theta$, and $Q$ is an integer such that
$Q=\frac{\cos \omega_0}{\eps}+O(\eps^{-\theta})$. Then one has
\begin{equation*}
\tilde{F}_{\eps,I,Q}(t)=c_I H\left( \frac{t}{\cos \omega_0+\sin
\omega_0} \right) +O\big( \EE_{\theta,\theta_1,\delta} (\eps)\big)
\end{equation*}
uniformly for $t$ in a compact subset of $(0,1)\cup (1,2)\cup
(2,\infty)$, with $\EE_{\theta,\theta_1,\delta}(\eps)$ as in
Proposition \ref{P3.1}.
\end{prop}

{\sl Proof.} We write $I=[\tan \omega_0,\tan \omega_1]$. Since
$\vert R_{1/Q}(\omega)-(q+a)\vert \leq 1$ whenever
$\eps_{\omega,Q}=(q,a)$, it follows that
\begin{equation*}
\vert R_h (\omega)-l_h(\omega)(\cos \omega+\sin \omega)\vert \leq
1,\qquad \forall h,\ \forall \omega .
\end{equation*}
This manifestly implies for all $\omega \in [\omega_0,\omega_1]$
\begin{equation}\label{3.12}
\begin{split}
 & \tilde{H}_{\eps,I,Q} \left( \frac{t+\eps}{\cos \omega_1+\sin
\omega_1} \right) \leq \tilde{H}_{\eps,I,Q} \left(
\frac{t+\eps}{\cos \omega +\sin \omega} \right) \\
& \qquad =\left| \left\{ \omega \, ;\, \tan \omega \in I,\
l_{\frac{1}{Q}} (\omega) (\cos \omega+\sin \omega)
>\frac{t}{\eps}+1 \right\} \right| \\
& \qquad \qquad \leq \tilde{F}_{\eps,I,Q}(t) \\
& \qquad \leq \left| \left\{ \omega \, ;\, \tan \omega \in I,\
l_{\frac{1}{Q}} (\omega) (\cos \omega+\sin
\omega)>\frac{t}{\eps}-1 \right\} \right| \\
& =\tilde{H}_{\eps,I,Q} \left( \frac{t-\eps}{\cos \omega+\sin
\omega} \right) \leq \tilde{H}_{\eps,I,Q} \left(
\frac{t-\eps}{\cos \omega_0+\sin \omega_0} \right) .
\end{split}
\end{equation}

By Proposition \ref{P3.1}, the mean value theorem and
\begin{equation*}
\frac{t+\eps}{\cos \omega_1+\sin \omega_1} -\frac{t}{\cos
\omega_0+\sin \omega_0} =O(\eps^\theta),
\end{equation*}
we get
\begin{equation}\label{3.13}
\begin{split}
& \tilde{H}_{\eps,I,Q} \left( \frac{t+\eps}{\cos \omega_1+\sin
\omega_1} \right)  =c_I H \left( \frac{t+\eps}{\cos \omega_1+\sin
\omega_1} \right)+O\big( \EE_{\theta,\theta_1,\delta}(\eps)\big)
\\ & \qquad =c_I \bigg( H\Big( \,\frac{t}{\cos \omega_0+\sin \omega_0}
\Big)+O(\| H\|_\infty \eps^\theta) \bigg) + O\big(
\EE_{\theta,\theta_1,\delta}(\eps)\big)
\\ & \qquad = c_I H\left( \frac{t}{\cos \omega_0+\sin \omega_0}\right) +
O\big( \EE_{\theta,\theta_1,\delta}(\eps)\big).
\end{split}
\end{equation}
In a similar way one gets
\begin{equation}\label{3.14}
\tilde{H}_{\eps,I,Q} \left( \frac{t-\eps}{\cos \omega_0+\sin
\omega_0} \right)  =c_I H\left( \frac{t}{\cos \omega_0+\sin
\omega_0}\right) + O\big( \EE_{\theta,\theta_1,\delta}(\eps)\big).
\end{equation}
The statement follows from \eqref{3.12}, \eqref{3.13} and
\eqref{3.14}. \qed

\bigskip

\section{Proof of Theorems \ref{T1.1} and \ref{T1.2}}
\begin{figure}[ht]
\begin{center}
\includegraphics*[scale=0.7, bb=0 0 240 220]{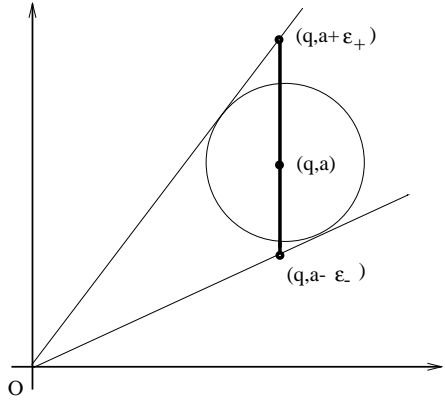}
\end{center}
\caption{A circular scatterer} \label{Figure7}
\end{figure}

We focus now on the case of circular scatterers of
radius $\eps$. For each integer lattice point
$(q,a)\in \Z^{2\ast}$, let $(q,a\pm \eps_\pm)$ be the intersection
of the line $x=q$ with the tangents from $O$ to the circle
of center $(q,a)$ and radius $\eps$ (see Figure \ref{Figure7}).

The quantities $\eps_\pm$ are computed from
\begin{equation*}
\eps=\frac{\big| a-\frac{a\pm \eps_\pm}{q}\cdot q\big|}{\sqrt{
1+\big(\frac{a\pm \eps_\pm}{q}\big)^2}} =\frac{\eps_\pm q}{\sqrt{
q^2+(a\pm \eps_\pm )^2}} \, ,
\end{equation*}
which eventually gives
\begin{equation}\label{4.1}
\eps_\pm (q,a) =\frac{\eps \sqrt{q^2(q^2+a^2)-a^2
\eps^2}}{q^2-\eps^2} \pm \frac{a\eps^2}{q^2-\eps^2}
=\frac{\eps}{\cos \arctan \frac{a}{q}}+O\left( \frac{\eps^2}{q}
\right) .
\end{equation}

Let $I=[\tan \omega_0,\tan \omega_1]\subseteq [0,1]$ be an
interval of size $\vert I\vert \asymp \eps^\theta$ with
$0<\theta <\frac{1}{2}$. It follows
from \eqref{4.1} that there exists $\eps_0 =\eps_0 (\theta) >0$
such that for all
$\eps <\eps_0$ and all $\frac{a}{q} \in \FI$ one has
\begin{equation}\label{4.2}
\frac{\eps-\eps^{3/2}}{\cos \omega_0} \leq \eps_\pm (q,a)
\leq \frac{\eps+\eps^{3/2}}{\cos \omega_1} \, .
\end{equation}
Since $\theta<\frac{1}{2}$ it also follows that
\begin{equation*}
\frac{\cos \omega_0}{\eps-\eps^{3/2}} - \frac{\cos
\omega_1}{\eps+\eps^{3/2}} \ll \eps^{\theta-1} .
\end{equation*}
Thus we find integers $Q^-$ and $Q^+$ such that
\begin{equation}\label{4.3}
Q^- \leq \frac{\cos \omega_1}{\eps+\eps^{3/2}} \leq
\frac{\cos \omega_0}{\eps-\eps^{3/2}} \leq Q^+ \quad
\mbox{\rm and} \quad Q^+-Q^-\ll \eps^{\theta-1} .
\end{equation}

By \eqref{4.2} and \eqref{4.3} it follows for $\eps <\eps_0$
and $\frac{a}{q} \in \FI$ that
\begin{equation*}
\frac{1}{Q^+} \leq \eps_\pm (q,a)\leq \frac{1}{Q^-}
\end{equation*}
and
\begin{equation}\label{4.4}
Q^\pm =\frac{\cos \omega_0}{\eps}+O(\eps^{\theta-1}) .
\end{equation}

Obvious monotonicity properties of
\begin{equation*}
H_{\eps,I}(t)=\bigg| \bigg\{ \omega \, ;\, \tan \omega \in I,\
\tau_\eps (\omega)>\frac{t}{\eps} \bigg\} \bigg|
\end{equation*}
and of $\tilde{H}_{\eps,I,Q}(t)$ imply that
\begin{equation}\label{4.5}
H_{\eps,I}(t)\leq \bigg| \bigg\{ \omega \, ;\, \tan \omega \in I,\
l_{\frac{1}{Q^+}} (\omega)
>\frac{t}{\eps}-2\eps \bigg\} \bigg|
=\tilde{H}_{\eps,I,Q^+}(t-2\eps^2)
\end{equation}
and
\begin{equation}\label{4.6}
H_{\eps,I}(t)=\tilde{H}_{\eps,I,Q^-}(t+2\eps^2).
\end{equation}

Due to \eqref{4.4} we may apply Proposition \ref{P3.1} to get
\begin{equation}\label{4.7}
\tilde{H}_{\eps,I,Q^\pm} (t\mp 2\eps) =c_I H(t\mp
2\eps)+O(\eps^{\min
(2\theta,\frac{1}{2}-2\theta_1-\delta,\theta+\theta_1-\delta)}).
\end{equation}

In view of \eqref{4.5}, \eqref{4.6} and \eqref{4.7} we gather for
any $t\neq 1,2$
\begin{equation}\label{4.8}
H_{\eps,I}(t)=c_I H(t)+O(\eps^{\min
(2\theta,\frac{1}{2}-2\theta_1-\delta,\theta+\theta_1-\delta)}) .
\end{equation}

Finally we choose a partition $[0,1]=\cup_{j=1}^N I_j$ with
intervals $I_j$ of equal size $\vert I_j \vert =\frac{1}{N} \asymp
\eps^\theta$. Summing in \eqref{4.8} over $I\in \{
I_1,\dots,I_N\}$ we arrive at
\begin{equation*}
H_\eps (t)=\frac{4H(t)}{\pi} \sum\limits_{j=1}^N c_j +O(\eps^{\min
(\theta,\frac{1}{2}-2\theta_1-\theta-\delta,\theta_1-\delta)}).
\end{equation*}

The proof of Theorem \ref{T1.1} is complete once we choose
$\theta=\theta_1=\frac{1}{8}$.

\medskip

{\sl Proof of Theorem \ref{T1.2}.} Arguing as above we get from
Proposition \ref{P3.7} taking $\theta=\theta_1=\frac{1}{8}$:
\begin{equation*}
F_{\eps,J}(t)=\tilde{F}_{\eps,J,Q^\pm} (t\mp 2\eps^2)
=c_J H\left( \frac{t}{\cos \omega_0+\sin \omega_0}\right)
+O(\eps^{\frac{1}{4}-\delta})
\end{equation*}
for any interval $J=[\tan \omega_0,\tan \omega_1] \subseteq
[0,1]$. Thus if the intervals
$I_j=[\tan \omega_j,\tan \omega_{j+1}]$, $1\leq
j\leq N=[\eps^{-\frac{1}{8}}]$, are such that $\vert I_j\vert
=\frac{1}{N} \asymp \eps^\frac{1}{8}$, then
\begin{equation*}
F_\eps (t) =\frac{4}{\pi} \sum\limits_{j=1}^N F_{\eps,I_j}(t)
=\frac{4}{\pi} \sum\limits_{j=1}^N (\omega_{j+1}-\omega_j) H
\left( \frac{t}{\cos \omega+\sin \omega} \right) +O
(\eps^{\frac{1}{8}-\delta}).
\end{equation*}

It only remains to compare the main terms in the relation above
and in Theorem \ref{T1.2}. This is
achieved by merely appying the mean value theorem:
\begin{equation*}
\begin{split}
\frac{4}{\pi} \int\limits_0^{\pi/4} H\left( \frac{t}{\cos
\omega+\sin \omega} \right) & =\frac{4}{\pi} \sum\limits_{j=1}^N
\int\limits_{\omega_j}^{\omega_{j+1}} H\left( \frac{t}{\cos \omega
+\sin \omega}\right) d\omega \\
& =\frac{4}{\pi} \sum\limits_{j=1}^N (\omega_{j+1}-\omega_j)
\bigg( H\Big( \frac{t}{\cos \omega_j+\sin \omega_j} \Big)
+O(\eps^\frac{1}{8})\bigg) \\
& =\frac{4}{\pi} \sum\limits_{j=1}^N (\omega_{j+1}-\omega_j)
H\left( \frac{t}{\cos \omega_j+\sin \omega_j} \right)
+O(\eps^\frac{1}{8}) ,
\end{split}
\end{equation*}
which concludes the proof of Theorem \ref{T1.2}. \qed

\bigskip

\bigskip

\end{document}